\documentclass{amsart}

\calclayout

\usepackage[unicode]{hyperref}
\usepackage[lite]{amsrefs}

\def\paragraph{\subsection}
\makeatletter
\def\@sect#1#2#3#4#5#6[#7]#8{%
  \edef\@toclevel{\ifnum#2=\@m 0\else\number#2\fi}%
  \ifnum #2>\c@secnumdepth \let\@secnumber\@empty
  \else \@xp\let\@xp\@secnumber\csname the#1\endcsname\fi
  \@tempskipa #5\relax
  \ifnum #2>\c@secnumdepth
    \let\@svsec\@empty
  \else
    \refstepcounter{#1}%
    \edef\@secnumpunct{%
      \ifdim\@tempskipa>\z@ % not a run-in section heading
        \@ifnotempty{#8}{.\@nx\enspace}%
      \else
        \@ifempty{#8}{.}{.\@nx\enspace}%
      \fi
    }%
      \ifnum #2=\tw@ \def\@secnumfont{\bfseries}\fi{}%
    \protected@edef\@svsec{%
      \ifnum#2<\@m
        \@ifundefined{#1name}{}{%
          \ignorespaces\csname #1name\endcsname\space
        }%
      \fi
      \@seccntformat{#1}%
    }%
  \fi
  \ifdim \@tempskipa>\z@ % then this is not a run-in section heading
    \begingroup #6\relax
    \@hangfrom{\hskip #3\relax\@svsec}{\interlinepenalty\@M #8\par}%
    \endgroup
    \ifnum#2>\@m \else \@tocwrite{#1}{#8}\fi
  \else
  \def\@svsechd{#6\hskip #3\@svsec
    \@ifnotempty{#8}{\ignorespaces#8\unskip
       \@addpunct.}%
    \ifnum#2>\@m \else \@tocwrite{#1}{#8}\fi
  }%
  \fi
  \global\@nobreaktrue
  \@xsect{#5}}
\makeatother

\makeatletter
\def\pxspace{\@ifnextchar.{\@}{.\@\xspace}}
\makeatother

\newcommand{\mvert}{\,|\,}      % "such that":  \{ x \in S \mvert x > y \}

\newcommand{\problem}[1]{\pagebreak[0]\hypertarget{problem-#1}{\textbf{#1.}}}

\begin{document}

\title[Tableaux open problems]%
      {On the combinatorics of tableaux --- A notebook of open problems}
\author{Dale R. Worley}
\email{worley@alum.mit.edu}
\date{\today} % format is Mmm dd, yyyy.

\begin{abstract}
Inspired by the the Kourovka Notebook of unsolved problems in group
theory\cite{KhukhMaz2024a}, this is a notebook of unsolved problems in
the combinatorics of tableaux.  Contributions to the notebook are
invited.
\end{abstract}

\maketitle

%\tableofcontents
%\listoffigures

\begin{quotation}
Mathematicians can be subdivided into two types: problem solvers
and theorizers. Most mathematicians are a mixture of the two although
it is easy to find extreme examples of both types.
---~Gian-Carlo Rota, "Problem Solvers and Theorizers" in \cite{Rot1997b}
\end{quotation}

% Push the introduction down a bit.
\vspace{2em}

\textit{Version 3 (5 April 2026)}
\vspace{\baselineskip}

% Be careful editing this paragraph because the long URL screws up
% TeX's line splitting algorithm.
Contributions to the notebook are invited.
Please send all contributions, comments, and corrections to the author.
Individual problems have HTML anchors of the form
``\texttt{nameddest=problem-}\textit{n}'', so you can reference an
individual problem with a URL like
``\texttt{https://.../problem-notebook.pdf\#nameddest=problem-10}''.

% https://alum.mit.edu/www/worley/Math/problem-notebook.pdf\#nameddest=problem-xxx

\section{Problems of 2024}

\problem{1} \cite{Stan1988a}*{sec.~6 Prob.~1} Classify all differential posets.
\cite{Fom1988}\linebreak[0]%
\cite{Stan2012a}*{sec.~3.21}\linebreak[0]%
\cite{Byrn2012a}*{Ch.~6}

\problem{2} \cite{Wor} Consider the tree $SSYT$ of shifted standard Young
tableaux, with
the tree structure generated by containment.  Computation suggests
that the only automorphism of $SSYT$ (as an unlabeled tree) is the
identity.  Computation further suggests that each node of $SSYT$ is
uniquely determined by the census (cardinalities) of the ranks of the
subtree rooted at the node.

\problem{3} \cite{Wor} Consider the tree $SYT$ of (unshifted) standard Young
tableaux, with
the tree structure generated by containment.  For any node with a
symmetric shape, there is an automorphism of $SYT$ (as an unlabeled
tree) that transposes all
nodes in the subtree rooted at the node and leaves all other nodes fixed.
Computation suggests
that all automorphisms of $SYT$ are compositions of these transpose
automorphisms (with the identity being the composition of zero
transpose automorphisms).
Computation further suggests that each node of $SYT$ is
determined up to transpose by the census (cardinalities) of the ranks of the
subtree rooted at the node.

\section{Problems of 2025}

\problem{4} \cite{Fom1994a}*{Exam.~2.2.4} notes that lattice of Young
diagrams with $\leq r$ rows is differential with degree $r$.  Thus
the general theory produces an RSK algorithm from $r$-colored
permutations into pairs of such tableaux with suitable weighting.
This should be expandable into similar theories as are developed for
Young diagrams and shifted Young diagrams.  In particular, (1) Can jeu de
taquin be defined for these diagrams?  (2) Can a Greene invariant be
defined for $r$-colored permutations that groups permutations that
produce the same $P$ tableau?

\problem{5} \cite{Wor} Consider the process of randomly choosing a
Young diagram of
size $n$.  If we let $n \rightarrow \infty$ and scale the diagram down
by $\sqrt{n}$ in both dimensions, the average of the diagrams converges
to a limit shape; points within the shape have probability 1 of being
in the random diagram and points outside the shape have probability 0.
Surprisingly, the same limit shape results under three distributions
on the set of diagrams:  uniform (weight $1$), uniform on Young
tableaux (weight $f_\lambda$, the number of tableaux of shape
$\lambda$), and the Plancherel measure (weight $f_\lambda^2$).

Perform the parallel analysis on shifted Young diagrams, the set of
partitions into distinct parts.  Since the RSK algorithm for shifted
diagrams has a weighting
($c_\lambda = 2^{|\lambda| - \ell(\lambda)}$),
there are now six plausible distributions:  $1$, $g_\lambda$,
$g_\lambda^2$, $c_\lambda$, $g_\lambda c_\lambda$, and
$g_\lambda^2 c_\lambda$, where $g_\lambda$ is the number of shifted
Young tableaux (without coloring) of shape $\lambda$.

\problem{6} \cite{Wor} Develop a similar theory of limit shapes for
Young diagrams
with $\leq r$ rows.  The initial challenge is studying the general
asymptotic behavior of the $r-row$ Young diagrams generated by
$r$-colored permutations so that a proper scaling can be defined so
that the required limit is defined.

{\newcommand{\mbbYF}{\mathbb{YF}}
\problem{7} \cite{Wor} Develop a similar theory of limit shapes for
Young--Fibonacci diagrams.  The initial challenge is defining a
concept of scaling between smaller and larger diagrams.  E.g., if the
scaling is a doubling, we need to define an embedding $D$ of $\mbbYF_r$
into $\mbbYF_r$ that maps words (diagrams) of size $k$ into words of
size $2k$, $\mbbYF_r^{(k)} \rightarrow \mbbYF_r^{(2k)}$.
The embedding should be ``nearly isometric'' in that it maps words that are
``near'' each other in $\mbbYF_r$ (in the distance metric of the
Hasse diagram) into words that are near each other $\mbbYF_r$, and
vice-versa.  And every word in $\mbbYF_r$ must be near an element of the
image of $D$.

Given such a doubling map, it defines a limit space within which the
question can be asked whether randomly-chosen $r$-colored permutations
have a limit shape.  If there is a limit shape, determine it.
Since $\mbbYF_r$ as dual graded graphs has a unique weighting 1,
there are no colorings of the tableaux, so there are only three
distributions of diagrams to be considered.
}

\problem{8} \cite{Wor} Consider the tree $k$-$SYT$ of standard Young tableaux
with at most $k$ rows, with the tree structure generated by
containment.  What are its automorphisms, when it is considered an
unlabeled tree?

\problem{9} \cite{Wor} Consider the tree $SYFT$ of standard ``Young--Fibonacci
tableaux'', with the tree structure generated by
containment, or equivalently, the tree of upward paths in the
Young--Fibonacci lattice.
What are the automorphisms of this tree, when it is considered an
unlabeled tree?

\problem{11} \cite{Stan2017a}*{sec.~1} via \cite{Set}
Define $\nu_w = \mathfrak{S}_w(1,\ldots,1)$ which is the sum of the
coefficients of the Schubert polynomial $\mathfrak{S}_w$ of the
permutation $w$.  How large can $\nu_n$ be for $w \in S_n$?
Define $u(n)$ to be this maximum.  We have some rather crude bounds
which show that,
\begin{equation*}
\frac{1}{4} \leq \liminf_{n \rightarrow \infty} \frac{\log_2 u(n)}{n^2}
\leq \limsup_{n \rightarrow \infty} \frac{\log_2 u(n)}{n^2} \leq \frac{1}{2}.
\end{equation*}
Determine better bounds.
That is, can we find $1/4 \leq \alpha \leq 1/2$ such that
$u(n) \sim 2^{\alpha n^2}$?

\cite{MoraPakPan2019a}*{Th.~1.3} improves the lower bound to
$\approx 0.2932362762$.  \cite{MoraPanPet2025a}*{Rem.~6.8} improves
the upper bound to 0.37.

\problem{12} \cite{Stan2017a}*{sec.~1}\cite{MerSmir2015a}*{Ques.~5.6}
via \cite{Set} Define
$\nu_w = \mathfrak{S}_w(1,\ldots,1)$ which is the sum of the coefficients of
the Schubert polynomial $\mathfrak{S}_w$ of the permutation $w$.  For
a given $n$, define $u(n)$ to be the maximum $\nu_w$ for $w \in S_n$.
For a given n, what permutations $w \in S_n$ maximize $\nu_w$?
As $n \rightarrow \infty$, is there some kind of ``limiting shape'' of
the permutations $w \in S_n$ for which $\nu_w = u(n)$?

\problem{13} \cite{MerSmir2015a}*{Conj.~5.7} via \cite{Set} Define
$\nu_w = \mathfrak{S}_w(1,\ldots,1)$ which is the sum of the
coefficients of the Schubert polynomial $\mathfrak{S}_w$ of the
permutation $w$.  For a given n, prove that the permutation $w \in S_n$
that maximizes $\nu_w$ is ``Richardson'' or ``layered'' in that it has
the form [Def.~5.5]
$(i_1, i_1-1, \ldots, 2, 1, i_2, i_2-1, \ldots, i_1+1, i_3, \ldots, i_2+1, \ldots)$
for some $i_1 < i_2 < i_3 < \cdots <i_k = n$.

\problem{14} \cite{MoraPakPan2019a}*{Conj.~4.1} via \cite{Set}
Define $\Upsilon_w = \mathfrak{S}_w(1,\ldots,1)$ which is the sum of the
coefficients of the Schubert polynomial $\mathfrak{S}_w$ of the
permutation $w$.  For permutations $w$ and $v$, define $w \otimes v$ to be
the permutation whose permutation matrix is the Kronecker product of
the permutation matrices of $w$ and $v$.  Define $1^n$ as the identity
permutation in $S_n$.  Is it true that for any integer
$k \geq 2$, $\Upsilon_{w \otimes 1^k} \geq \Upsilon_w^{k^2}$?

\problem{15} \cite{Byrn2012a}*{Conj.~6.6} Is it true that an
$r$-differential poset $P$ is a lattice if and only if $P$ does not
contain a crown covering a crown?  (This is known to be true for $r=1$.)

\problem{16} \cite{Byrn2012a}*{Conj.~6.7} What are all the
1-differential posets that do not contain [a particular induced
subposet]?

\problem{17} \cite{Stan1988a}*{sec.~6 Prob.~6}
Fix $r \in \mathbb{P}$. What is the least number
of elements of rank $n$ that an $r$-differential poset can have?
It seems plausible that the minimum value is achieved by
$\mathbb{Y}^{\times r}$, the $r$-th cartesian power of the lattice of partitions.
(The maximum value is attained by $\mathbb{YF}_r$, the
$r$-differential Young--Fibonacci lattice.\cite{Byrn2012a}*{Th.~1.2})

\problem{18} \cite{Stan2012a}*{Ch.~3 Exer.~29}
A finite lattice $L$ has $n$ join-irreducibles.  What is the most
number $f(n)$ of meet-irreducible elements $L$ can have?

\problem{19} \cite{Stan2012a}*{Ch.~3 Exer.~135(b)}
Let $\Lambda_n$ be the set of all $p(n)$ partitions of the integer
$n \geq 0$.  Order $\Lambda_n$ by refinement.  Determine the M\"obius
function $\mu(\lambda,\rho)$ of $\Lambda_n$.  (This is trivial when
$\lambda = \left<1^n\right>$ and easy when
$\lambda = \left<1^{n-2}2^1\right>$.)

\problem{20} \cite{Stan2012a}*{Ch.~3 Exer.~155(c)}
Find all finite modular lattices for which every interval is
self-dual.
Do the same for finitary modular lattices with finite covers.

\problem{21} \cite{Stan2012a}*{Ch.~3 Exer.~183(f)}
Let $\mathfrak{S}_n$ be the set of permutations of $n$ ordered by the
strong Bruhat order.  Characterize the intervals of $\mathfrak{S}_n$
that are boolean algebras and compute their total number.

\problem{22} \cite{Stan2012a}*{Ch.~3 Exer.~185(j)}
Let $\mathfrak{S}_n^W$ be the set of permutations of $n$ ordered by the
weak Bruhat order.  Characterize the intervals of $\mathfrak{S}_n^W$
that are distributive lattices and compute their total number.
The values for $1 \leq n \leq 8$ are
1, 2, 16, 124, 1262, 15898, 238572, 4152172.

\problem{23} \cite{Stan2012a}*{Ch.~3 Exer.~215(c)}
Let $P$ be an $r$-differential poset and let $p_i$ be the size of rank
$i$ of $P$.  Show that $p_i < p_{i+1}$ except for the case $i=0$ and
$r=1$.

\problem{24} \cite{Stan1999a}*{Ch.~6 Exer.~25(i)}
Let the symmetric group $\mathfrak{S}_n$ act on the polynomial ring
$A = \mathbb{C}[x_1,\ldots,x_n,y_1,\ldots,y_n]$ by permuting both the
$x_\bullet$ and the $y_\bullet$ simultaneously;
$w \cdot f(x_1,\ldots,x_n,y_1,\ldots,y_n) =
f(x_{w(1)}, \ldots, x_{w(n)},y_{w(1)}, \ldots, y_{w(n)})$ for all
$w \in \mathfrak{S}_n$.
Let $I$ be the ideal generated by all invariants of positive degree,
i.e.,
$$ I = \left< f \in A: w \cdot f = f \textup{ for all }
w\in\mathfrak{S}_n \textup{, and } f(0)=0 \right>. $$
Prove that the Catalan number $C_n = \frac{1}{n+1} \binom{2n}{n}$
is the dimension of the subspace of $A / I$ affording the sign
representation, i.e.,
$$ C_n = \operatorname{dim} \{ f\in A/I:
w \cdot f = (\operatorname{sgn} w) f \textup{ for all }
w \in \mathfrak{S}_n \}. $$

\problem{25} \cite{Stan1999a}*{Ch.~7 Exer.~33(c)}
Can anything be said in general about the number of distinct monomials
in the Schur function $s_\lambda(x_1,\ldots,x_n)$ for arbitrary
$\lambda$?

\problem{26} \cite{Stan1999a}*{Ch.~7 Exer.~38(b)}
Fix $0 \leq k \leq \binom{n}{2}$, and let $\ell(w)$ denote the number
of inversions of the permutation $w\in\mathfrak{S}_n$.  Let $\lambda$
and $\mu$ b partitions of length at most $n$, with
$\mu\subset\lambda$.  Define the symmetric function
$$ t_{\lambda/\mu,k} =
(-1)^k \sum_{\overset{w \in\mathfrak{S}_n}{\ell(w)\geq k}}
\varepsilon_w h_{\lambda+\delta-w(\mu+\delta)}. $$
Thus $t_{\lambda/\mu,k}$ is a ``truncation'' of the Jacobi--Trudi
expansion \cite{Stan1999a}*{Ch.~7 eq.~7.69} of $s_{\lambda/\mu}$.
Is there a ``nice'' combinatorial interpretation of the scalar product
$\left< t_{\lambda/\mu,k}, s_\nu \right>$?

\problem{27} \cite{Stan1999a}*{Ch.~7 Exer.~55(b)}
Let $\rho^\lambda: \mathfrak{S}_n \rightarrow
\operatorname{GL}(m,\mathbb{C})$ be an irreducible representation of
$\mathfrak{S}_n$ with character $\chi^\lambda$ (so $m = f^\lambda$).
Is it possible to count the number of $\lambda$'s for which
$\rho^\lambda(\mathfrak{S}_n) \subset \operatorname{SL}(m,\mathbb{C})$?

\problem{28} \cite{Stan1999a}*{Ch.~7 Exer.~68(b)}
Let $G$ be a finite group of order $g$.  Given $w \in G$, let $f(w)$ be the
number of pairs $(u,v) \in G \times G$ satisfying
$w = uvu^{-1}v^{-1}$ (the commutator of $u$ and $v$).
Thus $f$ is a class function on $G$ and hence a linear combination
$\sum c_\chi \chi$ of irreducible characters $\chi$ of $G$.
The multiplicity $c_\chi$ of $\chi$ in $f$ is equal to $g/\chi(1)$.
Since $\chi(1)| g$, it follows that $f$ is a character of $G$.
Find an explicit $G$-module $M$ whose character is $f$, which directly
proves that $f$ is a real character.
Since in the equation $f = \sum \frac{g}{\chi(1)} \chi$ all the
$\chi$'s are linearly independent, this would provide a
new proof of the basic result that $\chi(1)$ divides the order of $G$.

\problem{29} \cite{Stan1999a}*{Ch.~7 Exer.~71(e)}
Given a finite group $G$,
define $\psi_G$ as the character of the action of $G$ on itself by
conjugation.  That is, the character of permutation representation
$\rho:G \rightarrow \mathfrak{S}_G$ defined by
$\rho(x)(y) = x y x^{-1}$, where $\mathfrak{S}_G$ is the group of
permutations of $G$.
Define $\kappa_{n\lambda} = \left< \psi_{\mathfrak{S}_n}, \chi^\lambda \right>$,
where $\chi^\lambda$ is the irreducible representation of $\mathfrak{S}_n$
indexed by the partition $\lambda$.
$\kappa_{n\lambda} > 0$ with the sole exception of
$n=2, \lambda = (1,1)$.
Is there a ``nice'' combinatorial interpretation of the numbers
$\kappa_{n\lambda}$?

\problem{31} \cite{Bren2024a}*{Conj.~1.4} Lusztig's ``Combinatorial
Invariance Conjecture'' about the Kazhdan--Lusztig polynomials
$P_{x,y}(q)$:  Let $(W_1, S_1), (W_2, S_2)$ be two Coxeter systems,
and $u \leq v$ in $W_1$, $w \leq z$ in $W_2$ (in the strong Bruhat
order) be such that $[u, v] \simeq [w, z]$ (the intervals in the Strong
Bruhat order are isomorphic as posets). Then prove
$$ P_{u,v}(q) = P_{w,z}(q). $$

\problem{32} \cite{Bjor} via \cite{Bren2024a}*{Prob.~1.9}
Let $(W,S)$ be a Coxeter system and $u \leq v$ in $W$.
Is it true that then there are a Coxeter system $(W^\prime,S^\prime)$
and $w^\prime \in W^\prime$ such that
$$ P_{u,v}(q) = P_{e^\prime,w^\prime}(q) $$
where $e^\prime$ is the identity in $W^\prime$?

\problem{33} \cite{Bren2024a}*{Prob.~1.10}
Find a combinatorial interpretation for Kazhdan-Lusztig polynomials.

\problem{34} \cite{Haim1985a}*{sec.~1.0}
Is the class of linear lattices self-dual?

\problem{35} \cite{Jons1953b}*{sec.~5}
Is the class of linear lattices a variety?  (That is, is it definable by
a set of equations.)

\problem{36} \cite{Wor2024b}*{quest.~3.65}
Are there lattices that obey all orders of the Arguesian law but are
not linear?  Are there such lattices that are finite?

\problem{37} \cite{Wor2024b}*{quest.~3.74}
Is the class of lattices that are isomorphic to a sublattice of the
lattice of normal subgroups of some group self-dual?

\problem{38} \cite{Jons1953b}*{sec.~5}
Is the class of lattices that are isomorphic to a sublattice of the
lattice of normal subgroups of some group a variety?

\problem{39} \cite{Wor2024b}*{quest.~3.80}
Is the class of lattices that are isomorphic to a sublattice of the
lattice of subgroups of some Abelian group self-dual?

\problem{40} \cite{Jons1953b}*{sec.~5}
Is the class of lattices that are isomorphic to a sublattice of the
lattice of subgroups of some Abelian group a variety?

\problem{41} \cite{Wor2024b}*{quest.~3.92}
Given a characteristic $p$, is the class of lattices that are
isomorphic to a sublattice of the lattice of subspaces of some vector
space over a field of characteristic $p$ a variety?

\problem{42} \cite{Stan1990a}*{sec.~2}
For any vector $\mathbf{r}$ there is at most one
$\mathbf{r}$-differential \emph{distributive} lattice $L(\mathbf{r})$.
It is probably hopeless to determine for which vectors $\mathbf{r}$
$L(\mathbf{r})$ exists.  If the elements of the vector increase
``sufficiently fast'' then $L(\mathbf{r})$ exists.  Is it possible to
describe the necessary rate of growth?

\problem{43}
A central point of \cite{LiangSag2024a}*{Th.~3.1} is that for any
partition
$\lambda = (\lambda_1, \lambda_2, \ldots \lambda_{\ell(\lambda)})$ with
$\lambda_1 \geq \lambda_2 \geq \cdots \geq \lambda_{\ell(\lambda)} > 0$,
if we define
$\lambda + i = (\lambda_1+i, \lambda_2+i, \ldots, \lambda_{\ell(\lambda)}+i)$,
then the set of partitions $(\lambda + i)_{i \geq 0}$ has the property
that for every $n$, the set of intervals in the partition lattice
$\{ [\lambda+n, \lambda+(n+k)] \,|\, k \geq 0\}$ are all isomorphic.
This series of elements allows their Order Ideal
Lemma (Lem.~1.1) to be applied.

Similar series of elements can be found in other distributive Fomin
lattices, the lattice of strict partitions, the lattice of $k$-row
partitions, and the cylindric partitions.  This is surprising since
this property does not seem to be forced by the hypotheses of distributive
Fomin lattices.

Can this observation be extended to the Young--Fibonacci lattices?
That is, are there sequences of elements in $\mathbb{YF}_r$,
$(p_i)_{i \geq 0}$ for which for every $n$, the set of intervals
$\{ [p_n, p_{n+k}]\,|\,k \geq 0 \}$ are all isomorphic?
If there are, can they be used to prove an analog of their Th.~3.1?
If so, does this provide a clue to a ``semi-distributive'' property
of $\mathbb{YF}_r$?

\problem{44} \cite{Stan1990a}*{sec.~2}
For any vector $\mathbf{r}$ there is at most one
$\mathbf{r}$-differential \emph{distributive} lattice $L(\mathbf{r})$.
If we consider the class of $\mathbf{r}$-vectors that are arithmetic
progressions with
a given difference $\Delta$, then the class $L(\mathbf{r})$ of
differential lattices that are differential posets with such
$\mathbf{r}$-vectors is closed under cartesian products.
Conversely, if the cartesian product of two sequentially differential
posets each with at least two ranks is a sequentially differential
poset, then both the factors and
the product have $\mathbf{r}$-vectors that are arithmetic progressions
and they all have the same difference $\Delta$.
Applying
\cite{Stan2012a}*{Exer.~3.51(a) soln.} shows that $\Delta \leq 0$.  If
$\Delta = 0$, they are ordinary differential posets.  But if $\Delta <
0$, they must be finite.  Do there exist any finite distributive
differential posts with $\mathbf{r}$-vectors that are arithmetic
progressions other than the cartesian
exponentials of the 2- and 3-element chains?\cite{Stan1990a}*{Exam.~2.5}
Can the cases with $\Delta < 0$ be fully classified?

\problem{45} \cite{Fom1994a}*{sec.~2.3}.
Given a graded vector space $V = \bigoplus_{i=0}^{\infty} V_i$ and a
graded linear operator $D$ on it where $D(V_0) = \{0\}$ and
$D(V_{i+1}) \subset D(V_{i-1})$, is it always possible to define a
bilinear product on $V$ that turns it into a graded algebra with
identity for which $D$ is a derivation --- $D(ab) = D(a)b + aD(b)$?
It seems that one can always define such a product, but there is much
freedom and no clear guarantee the product has any nice properties.
So let us add that the product must be associative.
If this product exists, it allows the graded linear operator
$D$ to be used to construct a dual graded graph by the construction in the
citation.

{\newcommand{\mbbY}{\mathbb{Y}}\newcommand{\mbbZ}{\mathbb{Z}}
\problem{46} \cite{Gaetz2018a}*{Conj.~1.6}
Let $\mathfrak{G}: \{e\} = G_0 \subset G_1 \subset G_2 \subset \cdots$
be a tower of groups and let $r \in \mbbZ_{>0}$.
Define $P$ to be the splitting diagram of the representations of the
tower of groups, with weights that are the multiplicities of the
splittings.
The tower is \emph{$r$-dual} if $P$ is a self-dual
graded graph with differential degree $r$.
The tower is \emph{$r$-differential} if the tower is $r$-dual and the
weights of $P$ are all 1, it is a differential poset.
Show that:
(a) If $\mathfrak{G}$ is $r$-differential then
$P(\mathfrak{G}) \cong \mbbY^{\otimes r}$ and there exists an abelian
group $A$ of order $r$, not depending on $n$, such that
$G_n \cong A \wr S_n$ for all $n \geq 0$.
(b) If $\mathfrak{G}$ is $r$-dual then
$P(\mathfrak{G}) \cong (d_1\mbbY)\times \cdots \times(d_k \mbbY)$ for
some $d_1, \ldots, d_k$ and there exists a group $H$ of order $r$, not depending
on $n$, such that $G_n \cong H \wr S_n$ for all $n \geq 0$.
}

\section{Problems of 2026}

\problem{47} \cite{ChoiNamSOh2019a}*{sec.~1} The ``tableau switching''
process has a correspondence with jeu de
taquin.\cite{BenSottStroom1996a} There is a ``shifted tableau
switching'' process that similarly corresponds to shifted jeu de
taquin.  However there is an alternative ``modified shifting tableau
switching'' process with certain better properties.  ``It would be
very nice to develop a modified version of shifted jeu de taquin which
plays the same role as that of the shifted jeu de taquin in the
shifted switching process.''

\problem{48} \cite{FaigHerr1981a} via \cite{Far}
shows that a finite modular lattice $L$ with join-irreducible elements
$J$ is characterized by the partial ordering induced on $J$ and a
function $\rho: J \times J \rightarrow 2^J:
(a, b) \mapsto \{c \in J \mvert c \leq a \vee b\}$;
if two modular lattices have isomorphic posets $J$ and functions
$\rho$, they are isomorphic.  Can this be extended to a classification
by describing necessary and sufficient criteria for $\rho$ for it
to correspond to a modular lattice?

\problem{49} \cite{Wor2026c}
Prove that
a locally-finite distributive lattice has a lattice of ideals of prime
filters that contains at most three connected components (one
isomorphic to the lattice, perhaps one containing the empty ideal, and
perhaps one containing the ideal of all primes) iff
the lattice has finite width  (contains no infinite
antichains).

\problem{50} \cite{Stan2012a}*{Exer.~3.51(a) soln.}\cite{Wor2026e}
\cite{Stan2012a} shows that there exists exactly one distributive,
$d$-differential lattice for any positive integer $d$.  That analysis
easily generalizes to an analysis of weighted differential
distributive lattices, and when we constrain the differential degree
and the weights to be positive integers, it can be used to demonstrate
non-trivial constraints on possible lattices.  Can this technique be
extended to prove all of the ``strong'' hypotheses of the classification in
\cite{Wor2026b} to create a classification of all positive-weight
differential distributive lattices with $\hat{0}$?

\section{Solved Problems of 2025}

\problem{10} \cite{Stan2017a}*{sec.~4 and Conj.~4.1} via \cite{Set}
Define $\nu_w = \mathfrak{S}_w(1,\ldots,1)$ which is the sum of the
coefficients of the Schubert polynomial $\mathfrak{S}_w$ of the
permutation $w$.  It is well-known that $\nu_w = 1$ iff $w$ is
132-avoiding (or ``dominant'').
Conjecture 4.1.: We have $\nu_w = 2$ if and only if $w = a_1\cdots a_n$ has
exactly one subsequence $a_i a_j a_k$ in the pattern 132, i.e.,
$a_i < a_k < a_j$.

Proved by \cite{Weig2018a}*{Cor.~1.3}.

\problem{30} \cite{Stem} via \cite{Bren2024a}*{Prob.~1.1}
Let $(W, S)$ be a Coxeter system.
Let $T = \{ wsw^{-1}: s \in S, w\in W \}$ be its set of reflections
and for $w \in W$ let $\ell(w)$ be the length of $w$, the length of
the shortest product of generators that $= w$.
Is it true that $$ \sum_{t\in T} x^{\ell(t)} $$
is a rational generating function?

Via \cite{Biag}: Proved by \cite{DeMan1999a}*{Th.~4.1}.  A more
explicit formula for affine Coxeter systems is given by
\cite{BiagHohlSass2026a}*{Th.~1.2}.

% Push the e-mail address down a bit.
\vspace{3em}

\end{document}